\documentclass[fleqn]{mat01}
\usepackage{times,mathtimy,amssymb,latexsym}
\begin{document}

\setcounter{page}{271} \firstpage{271}


\def\ch{\mbox{\rm ch}}
\def\e{\mbox{\rm e}}
\def\Sq{\mbox{\rm Sq}}
\def\Vect{\mbox{\rm Vect}}

\def\Hom{\operatorname{Hom}}
\def\Im{\operatorname{Im~}}

\def\Q{\mbox{$\mathbb{Q}$}}
\def\Z{\mbox{$\mathbb{Z}$}}
\def\Pa{\mbox{$\mathbb{P}$}}

\def\rema{\trivlist\item[\hskip\labelsep{\it Remark.}]}
\def\remarr{\trivlist\item[\hskip\labelsep{\it Remarks.}]}

\newtheorem{theore}{Theorem}
\renewcommand\thetheore{\arabic{section}.\arabic{theore}}
\newtheorem{theor}[theore]{\bf Theorem}
\newtheorem{coro}[theore]{\rm COROLLARY}
\newtheorem{propo}[theore]{\rm PROPOSITION}
\newtheorem{lem}[theore]{Lemma}
\newtheorem{exam}[theore]{Example}

\renewcommand\theequation{\arabic{section}.\arabic{equation}}

\title{On the classification of complex vector bundles of stable
rank}

\markboth{Constantin B\v{a}nic\v{a} and Mihai Putinar}{Complex
vector bundles of stable rank}

\author{CONSTANTIN B\v{A}NIC\v{A}$^{*}$ and MIHAI PUTINAR$^{1}$}

\address{$^{1}$Mathematics Department, University of California, Santa
Barbara, CA~93106, USA\\
\noindent E-mail: mputinar@math.ucsb.edu}

\volume{116}

\mon{August}

\parts{3}

\pubyear{2006}

\Date{MS received 12 April 2006}

\begin{abstract}
One describes, using a detailed analysis of Atiyah--Hirzebruch
spectral sequence, the tuples of cohomology classes on a compact,
complex manifold, corresponding to the Chern classes of a complex
vector bundle of stable rank. This classification becomes more
effective on generalized flag manifolds, where the Lie algebra
formalism and concrete integrability conditions describe in
constructive terms the Chern classes of a vector
bundle.\renewcommand\thefootnote{} \footnote{$^{*}$Since
deceased.}
\end{abstract}


\maketitle

\section{Introduction}

Let $X$ denote a finite CW-complex of dimension $n$. For a natural
$r$, one denotes by $\Vect_{\rm top}^{r} (X)$ the isomorphism
classes of complex vector bundles on $X$, of rank $r$. It is
well-known that the map:
\begin{equation*}
\Vect_{\rm top}^{[n/2]}(X) \rightarrow \Vect_{\rm
top}^{r}(X),\quad E \mapsto E\oplus 1_{r-[n/2]}
\end{equation*}
is onto whenever $r \geq [n/2]$ (the integer part of $n/2$).
Moreover, for $r \geq n/2$, the same map is bijective (see for
instance \cite{13}).

On the other hand, Peterson \cite{17} has proved that, if $r \geq
n/2$ and $H^{2q}(X, \mathbb{Z})$ has no $(q - 1)!$-torsion for any
$q \geq 1$, then two vector bundles $E$ and $E'$ of rank $r$ on
$X$ are isomorphic if and only if they have the same Chern
classes.

An intriguing question is how to determine the range of the total
Chern class map:
\begin{equation*}
c\hbox{\rm :}\ \Vect_{\rm top}^{r}(X) \rightarrow \prod_{q}
H^{2q}(X, \mathbb{Z}).
\end{equation*}
According to previous observations, this would imply (under the
above torsion conditions) a classification of all rank $r$ complex
vector bundles on $X$, for stable rank $r \geq n/2$.

A few partial answers to this question are known. For instance, a
classical result of Wu asserts that any couple of cohomology
classes $(c_{1}, c_{2})\in H^{2}(X, \mathbb{Z}) \times H^{4}(X,
\mathbb{Z})$ coincides with $c(E)$ for a vector bundle $E$, on a
basis $X$ of dimension less than or equal to 4. For complex
projective spaces the range of $c$ was computed by Schwarzenberger
(see Appendix~1, pp.~22 of \cite{12}) and Thomas \cite{21}.
A~proof of Schwarzenberger--Thomas result is also given in the
paper by Switzer \cite{20}. If $n \leq 7$ and $H^{7}(X,
\mathbb{Z})$ has no 2-torsion, then a triple $(c_{1}, c_{2},
c_{3}), c_{i}\in H^{2i}(X, \mathbb{Z})$ coincides with $c(E)$ for
a rank-3 vector bundle $E$ on $X$, if and only if $c_{3} \equiv
c_{1}c_{2} + \Sq^{2}c_{2}$ in $H^{6}(X, \mathbb{Z}_{2})$ (see
\cite{3}). Moveover, for a compact complex manifold $X$ of complex
dimension 3, one has the bijection:
\begin{align*}
c\hbox{\rm :}\ \Vect_{\rm top}^{3}(X) &\rightarrow \{(c_{1},
c_{2}, c_{3})\hbox{\rm ;}\ c_{i} \in H^{2i}(X, \mathbb{Z})\quad
\hbox{for} \ \ i = 1, 2, 3\\[.2pc]
&\quad\ \ \hbox{and} \ \ c_{3} \equiv c_{1}c_{2} +
c_{1}(X)c_{2}(\bmod\ 2)\}.
\end{align*}
Here $c_{1}(X)$ is the first Chern class of the complex tangent
bundle of $X$.

In the present paper, we study the classification of stable-rank
bundles on some compact complex manifolds. First, we extend the
classification to manifolds of dimensions 4 and 5 (see
Propositions~3.1 and 3.2). We then consider some particular
classes of complex manifolds of arbitrary dimensions, including
the rational homogeneous manifolds. We list below the results
obtained in this case (see Propositions~4.1--4.3 and
Theorems~5.1--5.3 for precise statements).

Let $X$ be a compact complex manifold of complex dimension $n$.
\begin{enumerate}
\renewcommand\labelenumi{(\alph{enumi})}
\leftskip .05pc

\item If $H^{*}(X, \mathbb{Z})$ has no torsion, then $\Vect_{\rm
top}^{n}(X)$ is a bijection with those $n$-tuples of cohomology
classes $(c_{1}, \dots, c_{n})$ satisfying:
\begin{equation*}
\hskip -1.25pc \int_{X}\ch (c_{1}, \dots, c_{n})\ch (\xi)td(X) \in
\mathbb{Z}
\end{equation*}
for any $\xi \in K(X)$.

\item Assume that $H^{*}(X, \mathbb{Z})$ has no torsion and that
$H^{2}(X, \mathbb{Z})$ generates multiplicatively $H^{\rm even}(X,
\mathbb{Z})$. Then $\Vect_{\rm top}^{n}(X)$ is in bijection with
the classes $(c_{1}, \dots, c_{n})$ satisfying
\begin{equation*}
\hskip -1.25pc \int_{X}\ch (c_{1}, \dots, c_{n})\e^{\xi}td(X) \in
\mathbb{Z},
\end{equation*}
for any $\xi \in H^{2}(X, \mathbb{Z})$.

\item Assume that $X = G/B$, where $G$ is a complex algebraic,
semi-simple, simply connected Lie group and $B$ is a Borel
subgroup. Then, $\Vect_{\rm top}^{n}(X)$ is isomorphic to the set
of all classes $(c_{1}, \dots, c_{n})$ with the property
\begin{equation*}
\hskip -1.25pc \int_{X}\ch(c_{1},\dots, c_{n})\ch({\mathcal
O}_{X_{w}})\e^{-\rho}\in\mathbb{Z},
\end{equation*}
for $w \in W$ (the Weyl group of $G$). Moreover, these integrality
conditions are equivalent to
\begin{equation*}
\hskip -1.25pc \int_{X}\ch(c_{1}, \dots, c_{n})\e^{\chi}\in
\mathbb{Z},
\end{equation*}
for any character $\chi$ of a fixed maximal torus contained in
$B$.

Here ${\mathcal O}_{X_{w}}$ denotes the structure sheaf of the
closed Schubert cell $X_{w}$ associated with $w\in W$. Also,
$\rho$ stands for the semisum of all positive roots of $B$.

\item If $X$ is as before and $G$ is a product of simple groups of
type SL or Sp, then the conditions are
\begin{equation*}
\hskip -1.25pc \int_{X}\ch (c_{1},\dots, c_{n})\e^{\chi - \rho}
\in \mathbb{Z},
\end{equation*}
where $\chi = \omega_{i_{1}}+\cdots + \omega_{i_{k}}, 1\leq i_{1}
\leq \cdots\leq i_{k} \leq r, 0 \leq k \leq n - 3$ and
$\omega_{1},\dots, \omega_{r}$ are the fundamental weights of $G$
(we have put by convention that $\chi = 0$ for\break $k = 0$).

\parindent 1pc Actually the classes $\ch ({\mathcal O}_{X_{w}})$ appearing in assertion
(c) can be related to the natural basis $([X_{w}])_{w \in W}$ of
$H^{*}(X, \mathbb{Z})$, by using Demazure's results on the
desingularization of Schubert cells \cite{8}. However we ignore
whether there exists an accessible transformation which relates
the two bases $(\ch({\mathcal O}_{X_{w}}))_{w\in W}$ and
$([X_{w}])_{w\in W}$ of $H^{*}(X, \mathbb{Q})$. This is discussed
in the final part of the paper.

\parindent 1pc The main tools in the sequel are the Atiyah--Hirzebruch spectral
sequence, the integrality condition derived from the
Atiyah--Singer index theorem and of course the specific properties
of flag manifolds. The present article is motivated by the still
unsolved problem concerning the existence of analytic structures
on complex vector bundles, and the earlier attempts to find a
solution in low dimension \cite{3} via a topological
classification of vector bundles at stable rank.\vspace{-.5pc}
\end{enumerate}

\section{The Atiyah--Hirzebruch spectral sequence and a technical
result}

\setcounter{equation}{0}

This preliminary section contains some consequences of the
Atiyah--Hirzebruch spectral sequence for topological $K$-theory.
Among them, the subsequent technical Lemma~2.1 is the basis of all
results contained in this paper.

Let $X$ be a finite CW-complex and let $X^{i}$ stand for its
$i$-dimensional skeleton. Corresponding to the natural filtration,
\begin{equation*}
X = X^{n}\supset X^{n-1} \supset \cdots \supset X^{0} \supset
X^{-1} = \varnothing,
\end{equation*}
there exists the Atiyah--Hirzebruch spectral sequence which
converges to $K(X)$. The second level terms of this sequence are
\begin{equation*}
E_{2}^{p, q} = H^{p} (X, K^{q}(pt.)) =
\begin{cases}
H^{p}(X, \mathbb{Z}), &q \ \ \hbox{even}\\[.2pc]
0, &q \ \ \hbox{odd}.
\end{cases}
\end{equation*}
The coboundaries of this spectral sequence are denoted by
$d_{r}^{p,q}$. Due to the concrete form of $E_{2}^{p,q}$ one
easily finds $d_{r}^{p,q} = 0$ for $r$ even. Moreover, all
coboundaries $d_{r}^{p,q}$ have torsion (see \cite{1} or
\cite{10}).

It was also proved by Atiyah and Hirzebruch \cite{1} that
$d_{r}^{p,q} = 0$ for all $p, q, r$ $(r\geq 2)$, whenever
$H^{*}(X, \mathbb{Z})$ has no torsion. Indeed, because
$d_{3}\hbox{\rm :}\ H(X, \mathbb{Z}) \rightarrow H^{+3}(X,
\mathbb{Z})$ has torsion, it follows that $d_{3} = 0$, and so on.

For later use, the following detailed analysis of the coboundaries
$d_{r}^{p, q}$ is needed. First, Atiyah and Hirzebruch \cite{1}
noticed that, for a given class $\alpha \in H^{2p}(X,
\mathbb{Z})$, if $d_{2k + 1}\alpha = 0$ for all $k \geq 1$, then
there exists a class $\xi \in K(X)$ such that
\begin{equation*}
\xi|X^{2p-1} \ \ \hbox{is trival and} \ \ c_{p}(\xi) = (p -
1)!\alpha
\end{equation*}
(see also ch.~VIII of \cite{10}). Since any $\xi \in K(X)$ can be
represented as $\xi = [E] - [1_{m}]$, it follows that there exists
a vector bundle $E$ on $X$ which is trivial on $X^{2p - 1}$ and
its $p$th Chern class is prescribed as $c_{p}(E) = (p -
1)!\alpha$.

The vanishing conditions $d_{2k + 1}\alpha = 0, k \geq 1$, are
superfluous in the case of a \hbox{CW-complex} $X$ with
torsion-free cohomology.

A deeper analysis of the cohomology operations $d_{r}^{p, q}$
leads to the next result of Buhstaber \cite{7}. Let $p, q$ be
positive integers such that $\dim X < p + q$ and let us denote
\begin{equation*}
m = \prod_{\substack{s \geq 2\\[.2pc] s {\rm prime}}} s^{\big[\frac{q - 1}{2s -
1}\big]}.
\end{equation*}
Then for every $\xi \in E_{2}^{p, 0}\cong H^{p}(X, \mathbb{Z})$,
the multiple $m\xi$ is annihilated by all coboundaries $d_{2k +
1}, k \geq 1$.

In particular, if $\dim X = 10$, then $12d_{2k + 1}^{4,\cdot} = 0$
for all $k \geq 1$.

Next we prove that under some natural nontorsion conditions, the
range of the Chern class map can be computed, at least in
principle. Throughout this paper, $\mathbb{Z}_{m}$ denotes the
finite cyclic group of order $m$ and $\dot{p}$ stands for the
class of the integer $p$ in $\mathbb{Z}_{m}$, or an operation
deduced from it.

\begin{lem}
Let $X$ be a topological space having the homotopy type of a
finite CW-complex of dimension $n$. Assume that the group
$H^{2q}(X, \mathbb{Z})$ has no $(q - 1)!$-torsion for $q\leq (n/2)
- 1$ and that $d_{2k + 1}^{2p, *} = 0$ for $k \geq 1$ and $p \geq
2$. For every $p \geq 1$ one denotes
\begin{align*}
C_{p} &= \{(c_{1},\dots, c_{p});\quad c_{i}\in H^{2i}(X,
\mathbb{Z}) \ \
\hbox{and there exists a vector bundle}\\[.2pc]
&\quad\,E \ \ \hbox{on} \ \ X \ \ \hbox{with} \ \ c_{i}(E) =
c_{i}, \quad 1 \leq i \leq p\}.
\end{align*}
Then
\begin{enumerate}
\renewcommand\labelenumi{\rm (\arabic{enumi})}
\leftskip .1pc

\item $C_{1} = H^{2}(X, \mathbb{Z}), C_{2} = H^{2}(X, \mathbb{Z})
\times H^{4}(X, \mathbb{Z})${\rm ,}

\item There are functions $\varphi_{p}\hbox{\rm :}\
C_{p}\rightarrow H^{2p +2} (X, \mathbb{Z}_{p!})$ with the
property{\rm :}
\begin{equation*}
\hskip -1.25pc C_{p+1} = \{(c_{1}, \dots, c_{p + 1});
(c_{1},\dots, c_{p}) \in C_{p}, \dot{c}_{p + 1} =
\varphi_{p}(c_{1},\dots,c_{p})\}
\end{equation*}
for every $p \geq 1$.\vspace{-.7pc}
\end{enumerate}
\end{lem}

\begin{proof}$\left.\right.$

\begin{enumerate}
\renewcommand\labelenumi{(\arabic{enumi})}
\leftskip .1pc

\item The first equality $C_{1} = H^{2}(X, \mathbb{Z})$ is
obvious. In order to prove the second equality, let $(c_{1},
c_{2})$ be an arbitrary element of $H^{2}(X, \mathbb{Z}) \times
H^{4}(X, \mathbb{Z})$. There exists a complex line bundle $L$ on
$X$ with $c_{1}(L) = c_{1}$, and, according to the
Atiyah--Hirzebruch result mentioned above, there exists a vector
bundle $E$ on $X$ such that $c_{1}(E) = 0$ and $c_{2}(E) = c_{2}$.
Here we point out that the vanishing assumptions $d_{2k + 1}^{2p,
*} = 0$ are essential. Then the vector bundle $L\oplus E$ has the
first two Chern classes $c_{1}$ and $c_{2}$.

\item The function $\varphi_{p}\hbox{\rm :}\ C_{p}\rightarrow H^{2p+2}
(X, \mathbb{Z}_{p!})$ is defined as follows: Let $(c_{1},\dots,
c_{p})\in C_{p}$ and let $E$ be a vector bundle with $c_{i}(E) =
c_{i}, 1 \leq i \leq p$. Then
\begin{equation*}
\hskip -1.25pc \varphi_{p}(c_{1},\dots,c_{p}) =
\overline{c_{p+1}(E)},
\end{equation*}
that is the image of $c_{p+1}(E)$ through the natural map
$H^{2p+2}(X, \mathbb{Z})\rightarrow H^{2p+2}(X,
\mathbb{Z}_{p!})$.\vspace{-.5pc}
\end{enumerate}

In order to prove that this definition is correct, let $E'$ be
another vector bundle with $c_{i}(E') = c_{i}$ for $1 \leq i \leq
p$. We may assume in addition that $X$ is a CW-complex of
dimension $n$ and that $\hbox{rank}(E') = \hbox{rank}(E)$. Let
$X^{k}$ denote the $k$-skeleton of $X$. In view of Peterson's
theorem (see the introduction in \cite{17}), $[E| X^{2p+1}] = [E'|
X^{2p +1}]$, since $H^{s} (X, \mathbb{Z}) \cong H^{s}(X^{2p+1},
\mathbb{Z})$ for $s \leq 2p$ and $H^{s}(X^{2p+ 1}, \mathbb{Z}) =
0$ for $s > 2p + 1$. Consequently, the $K$-theory class $\xi = [E|
X^{2p+2}] - [E'| X^{2p +2}]\in \tilde{K} (X^{2p + 2})$ belongs to
the kernal of the restriction map $\rho$:
\begin{equation*}
\tilde{K} (X^{2p + 2}, X^{2p + 1})\xrightarrow{i} \tilde{K}(X^{2p
+ 2})\xrightarrow{\quad\rho\quad} \tilde{K}(X^{2p +1}).
\end{equation*}
From the above exact sequence of reduced $K$-theory one finds that
$\chi = i\eta$ for an element $\eta \in \tilde{K}(X^{2p + 2},
X^{2p+1})$.

But $X^{2p+2}/X^{2p+1}$ is a bouquet of $(2p + 2)$-dimensional
spheres, whence, by a theorem of Bott (see for instance
\cite{10}), $c_{i}(\eta) = 0$ for $1\leq i \leq p$ and $c_{p +
1}(\eta)$ is a multiple of $p!$. This implies the divisibility of
$c_{p + 1}(E|X^{2p+2}) - c_{p + 1}(E'| X^{2p+2})$ by $p!$. It
remains to remark that the restriction map $H^{2p+2}(X,
\mathbb{Z}_{p!})\longrightarrow H^{2p+2}(X^{2p+2},
\mathbb{Z}_{p!})$ is one-to-one, and this finishes the proof of
the correctness of the definition of $\varphi_{p}$.

Next we prove the equality in assertion (2) for $p \geq 1$.
Obviously $\varphi_{1} = 0$ and this identity agrees with $C_{2} =
H^{2}(X, \mathbb{Z}) \times H^{4}(X, \mathbb{Z})$. For an
arbitrary $p$, one remarks that the inclusion `$\subset$' in (2)
is a consequence of the definition of $\varphi_{p}$.

In order to prove the converse inclusion, let
$(c_{1},\dots,c_{p+1})$ be a $(p +1)$-tuple of cohomology classes
which satisfy: $(c_{1},\dots,c_{p}) \in C_{p}$ and $\dot{c}_{p+1}
= \varphi_{p}(c_{1},\dots,c_{p})$. Let $F$ be a vector bundle on
$X$ with the first Chern classes $c_{i}(F) = c_{i}, 1 \leq i \leq
p$. By the very definition of the function $\varphi_{p}$ one gets
$\dot{c}_{p+1} = \overline{c_{p_+1}(F)} =
\varphi_{p}(c_{1},\dots,c_{p})$. Hence there exists a cohomology
class $\sigma \in H^{2p+2}(X, \mathbb{Z})$, such that
\begin{equation*}
c_{p+1} - c_{p+1}(F) = (p!)\sigma.
\end{equation*}

The vanishing hypotheses in the statement assures the existence of
a vector bundle $G$ on $X$ with $c_{i}(G) = 0$ for $1\leq i \leq
p$ and $c_{p+1} (G) = (p!)\sigma$.

In conclusion, the vector bundle $E = F\oplus G$ has Chern classes
$c_{i}(E) = c_{i}$, $1 \leq i \leq p + 1$, and the proof of
Lemma~2.1 is over.\hfill $\Box$
\end{proof}

Let us remark that the sets $C_{p}$ and the functions
$\varphi_{p}$ are obviously compatible with morphisms $f\hbox{\rm
:}\ Y \rightarrow X$.

Lemma~2.1 will be used in the sequel in the following form: Let
$X$ and $\varphi_{1}, \varphi_{2},\dots$ be as before. Then the
cohomology clasess $(c_{1},\dots,c_{[n/2]})$ are the Chern classes
of a complex vector bundle on $X$ if and only if
\begin{equation*}
\dot{c}_{p+1} = \varphi_{p}(c_{1},\dots,c_{p}) \ \ \hbox{in} \ \
H^{2p+2}(X, \mathbb{Z}_{p!})
\end{equation*}
for every $p \geq 2$.

Indeed, if $\dot{c}_{3} = \varphi_{2}(c_{1}, c_{2})$, then
$(c_{1}, c_{2},c_{3})\in C_{3}$, so that $\varphi_{3}(c_{1},
c_{2},c_{3})$ makes a perfect sense, and so on. Of course, the
above vector bundle can be assume to be exactly of rank equal to
$[n/2]$.

The preceding result shows (at least theoretically) what kind of
congruences are needed to determine the range of the total Chern
class map. The main problem investigated in the following sections
is to establish concrete cohomological expressions for the
functions $\varphi_{p}$, for particular choices of the base space
$X$.

The simple applications of Lemma~2.1 appeared in \cite{3}. For the
sake of completeness, we include them in the present section,
too.\pagebreak

\begin{coro}$\left.\right.$\vspace{.5pc}

\noindent Let $X$ be a finite CW-complex of dimension less than or
equal to $7$ such that $H^{7}(X, \mathbb{Z})$ has no $2$-torsion.
Then $(c_{1}, c_{2}, c_{3})$ are the Chern classes of a rank $3$
vector bundle on $X$ if and only if
\begin{equation*}
c_{3}\equiv c_{1}c_{2} + \Sq^{2}c_{2}\quad \hbox{in} \ \ H^{6} (X,
\mathbb{Z}_{2}).
\end{equation*}
\end{coro}

\begin{proof}
The torsion conditions in Lemma~2.1 are automatically satisfied.
Moreover, the only possible non-zero coboundary in the
Atiyah--Hirzebruch spectral sequence is $d_{3}^{4}\hbox{\rm :}\
H^{4}(X, \mathbb{Z})\rightarrow H^{7}(X, \mathbb{Z})$.

But it is known that $d_{3}$ is a cohomological operation, which
coincides with Streenrod's square $\Sq^{3}$ (see \cite{10}).
Therefore $2\cdot d_{3}^{4} = 0$, which, in conjunction with the
lack of 2-torsion in $H^{7}(X, \mathbb{Z})$ shows that $d_{3}^{4}
= 0$.

At this moment we are able to apply Lemma~2.1 and we are seeking a
concrete form of the function $\varphi_{2}$. This is obtained by
Wu's formula:
\begin{equation}
w_{6} = w_{2}w_{4} + \Sq^{2}w_{4},
\end{equation}
referring to the universal Stiefel--Whitney classes (see for
instance, ch.~III.1.16 of \cite{13}).

Accordingly, the Chern classes $(c_{1}, c_{2}, c_{3})$ of a rank 3
vector bundle on $X$ satisfy
\begin{equation}
c_{3} \equiv c_{1}c_{2} + \Sq^{2}c_{2}\quad \hbox{in} \  \
H^{6}(X, \mathbb{Z}_{2}).
\end{equation}

In conclusion, the function $\varphi_{2}\hbox{\rm :}\ H^{2}(X,
\mathbb{Z}) \times H^{4}(X, \mathbb{Z}) \rightarrow H^{6}(X,
\mathbb{Z}_2)$ is precisely
\begin{equation*}
\varphi_{2}(c_{1}, c_{2}) = (c_{1}, c_{2} + \Sq^{2}c_{2})^{*}.
\end{equation*}

$\left.\right.$\vspace{-1.6pc}

\hfill $\Box$
\end{proof}

\begin{coro}$\left.\right.$\vspace{.5pc}

\noindent Let $X$ be a connected{\rm ,} compact complex manifold
of dimension $3$. Then
\begin{align*}
\Vect_{\rm top}^{3} (X)&\cong \{(c_{1}, c_{2}, c_{3}); c_{i} \in
H^{2i} (X,\mathbb{Z}), c_{1}c_{2} + c_{1}(X)c_{2}\equiv c_{3}
(\bmod\ 2),\\[.2pc]
&\,\quad i = 1, 2, 3\}.
\end{align*}
\end{coro}

Indeed{\rm ,} in this case $\Sq^{2}c_{2} = c_{1}(X)c_{2}(\bmod\
2)${\rm ,} where $c_{1}(X)$ stands for the first Chern class of
the tangent complex bundle of $X$.

The congruence relation $c_{3} \equiv c_{1} c_{2} + c_{1}(X)c_{2}
(\bmod\ 2)$ can alternatively be obtained in Corollary~$2.3$ from
the topological Riemann--Roch theorem.


For higher dimensional CW-complexes (i.e., of dimension greater
than 7) we ignore whether a finite result like Corollary~2.2
holds. It would be interesting for instance to identify the
functions $\varphi_3, \varphi_4$ among the identities satisfied by
the universal Chern classes and the secondary and tertiary
cohomological operations of order 6 and respectively 24.

For the rest of the paper we confine ourselves to only the simpler
case of compact complex manifolds.

\section{The stable classification on low dimensional compact
complex manifolds}

\setcounter{equation}{0}

\setcounter{theore}{0}

The Atiyah--Singer index theorem \cite{2} provides some
integrality conditions from which one can read the functions
$\varphi_3$ and $\varphi_4$, on compact complex manifolds of
dimensions 4 and 5 satisfy certain torsion-free cohomological
conditions.\pagebreak

\begin{propo}$\left.\right.$\vspace{.5pc}

\noindent Let $X$ be a compact{\rm ,} connected complex manifold
of complex dimension $4$ such that $H^6(X,\Z)$ and $H^7(X,\Z)$
have no $2$-torsion. Then $(c_1, c_2, c_3, c_4)$ are the Chern
classes of a {\rm (}unique{\rm )} rank $4$ vector bundle on
$X${\rm ,} if and only if{\rm :}
\begin{align}
c_3 &\equiv c_1 c_2 + \Sq^2 c_2\quad \hbox{in} \ \ H^6(X,
\Z_2),\\[.2pc]
c_4 &\equiv (c_1 + c_1(X))(c_3 - c_1c_2) + \frac{1}{2} (c_2(c_2 -
c_2(X))\nonumber\\[.2pc]
&\quad\, + c_1(X) (c_3 - (c_1 + c_1(X))c_2)) (\bmod 6).
\end{align}
\end{propo}

\begin{proof}
We will see below that (3.1) implies that the right-hand term of
(3.2) is an integer cohomology class in spite of the denominator
2.

The general torsion and vanishing assumptions in Lemma~1 are
fulfilled because $2d_3=0$ and $d_5^4=0$.

Let $E$ be a rank-4 complex vector bundle on $X$ with Chern
classes $c_i = c_i(E), i=1,2,3,4$. As before, formula (2.1) of Wu
implies the congruence (3.1), whence
\begin{equation*}
\varphi_2(c_1, c_2) = (c_1c_2 + \Sq^2 c_2)\quad \hbox{in} \ \
H^6(X, \Z_2).
\end{equation*}
The next function $\varphi_3$ is obtained from the Atiyah--Singer
index theorem.

First, recall the notation $\int_X \xi$ for the coefficient of the
component of degree 8 in the cohomology class $\xi \in H^* (X,
\Q)$, via the isomorphism $H^8(X, \Q) \cong \Q$ given by the
volume element.

The Atiyah--Singer index theorem \cite{2} implies that
\begin{equation}
\int_{X} \ch(E) td(X) \in \Z,
\end{equation}
for any vector bundle $E$ on $X$.

By denoting $\tau_i = c_i(T_X)$, where $T_X$ is the complex
tangent bundle of $X$, a straightforward computation shows that
relation (3.3) becomes
\begin{align}
&- \frac{1}{180} (\tau_1^4 - 4\tau_1^2 \tau_2 - 3\tau_2^2 - \tau_1
\tau_3 + \tau_4)\nonumber\\[.4pc]
&\quad\, + \frac{1}{24} c_1\tau_1\tau_2 + \frac{1}{24} (c_1^2 - 2
c_2)(\tau_1^2 + \tau_2)\nonumber\\[.4pc]
&\quad\, + \frac{1}{12} (c_1^3 - 3 c_1c_2 + 3c_3) + \frac{1}{24}
(c_1^4 - 4c_1^2 c_2 + 4c_1c_3 + 2 c_2^2 - 4c_4) \in \Z.
\end{align}

The last formula can be simplified by replacing $E$ with only a
few particular vector bundles. For instance, if $E=1_4$ (the
trivial bundle), then one finds that the first bracket is
divisible by 180. Then for $E=1_3\oplus L$, where $L$ is a line
bundle with $c_1(L)= c_1$, one gets
\begin{equation*}
\frac{1}{24} c_1 \tau_1\tau_2 + \frac{1}{24} c_2^2 (\tau_1^2 +
\tau_2) + \frac{1}{12} c_1^3 \tau_1 + \frac{1}{24} c_1^4 \in \Z.
\end{equation*}
By simplifying~(3.4) correspondingly one obtains exactly the
relation~(3.2).


Moreover, if $(c_{1}, c_{2}, c_{3})$ satisfies~(3.1), then by
Lemma~2.1, there exists a vector bundle $F$ on $X$ with $c_{i}(F)
= c_{i}, i = 1, 2, 3$. In view of the preceding computation one
finds the identity
\begin{align*}
&c_{4}(F) - (c_{1} + c_{1}(X))(c_{3} - c_{1}c_{2}) + 6\sigma\\[.2pc]
&\quad\,= \frac{1}{2}(c_{2}(c_{2} - c_{2}(X)) + c_{1}(X)(c_{3} -
(c_{1} + c_{1}(X))c_{2})),
\end{align*}
for some $\sigma \in H^{8}(X, \mathbb{Z})$. This shows that the
expression on the right-hand term always represents an integer
class of cohomology and that one can identify the function
$\varphi_{3}\hbox{\rm :}\ C_{3}\rightarrow H^{8}(X,
\mathbb{Z}_{6})$ with
\begin{align*}
\varphi_{3}(c_{1}, c_{2}, c_{3}) &= \bigg[ (c_{1} +
c_{1}(X))(c_{3} - c_{1}c_{2})\\[.4pc]
&\quad\,+\frac{1}{2}(c_{2}(c_{2} - c_{2}(X)) + c_{1}(X)(c_{3} -
(c_{1} + c_{1}(X))c_{2}))\bigg].
\end{align*}
Thus the proof of Proposition~3.1 is over.\hfill$\Box$
\end{proof}

By taking into account the properties of the Steenrod operation
Sq,
\begin{align*}
&\Sq^{2}\eta = c_{1}(X)\eta, \quad \eta\in H^{6}(X,
\mathbb{Z}_{2}),\\[.2pc]
&\Sq^{1}\xi = 0, \quad\xi\in H^{2}(X, \mathbb{Z}_{2}),\\[.2pc]
&\Sq^{2}(\xi c_{2}) = \xi \Sq^{2}c_{2} + \Sq^{1}\xi\cdot
\Sq^{1}c_{2} + \xi^{2}c_{2},
\end{align*}
one deduces from~(3.1) the relation \setcounter{equation}{4}
\begin{equation}
\xi(c_{1}c_{2} + c_{3} - (c_{1}(X) + \xi)c_{2}) \equiv 0\ (\bmod\
2)
\end{equation}
for every $\xi\in H^{2}(X, \mathbb{Z})$.

But $H^{6}(X, \mathbb{Z})$ was supposed to be without 2-torsion.
Hence $H^{3}(X, \mathbb{Z})$ has no 2-torsion, too, and the
restriction of scalars map $H^{2}(X, \mathbb{Z})\rightarrow
H^{2}(X, \mathbb{Z}_{2})$ is onto.

That means that relation~(3.5) is equivalent to~(3.1), because one
does not lose information by evaluating~(3.1) on the classes of
$H^{2}(X, \mathbb{Z}_{2})$.\hfill $\Box$\vspace{.7pc}

In conclusion, relations~(3.1) and (3.5) may be interchanged in
Proposition~3.1. Moreover, the above proof shows that
relation~(3.5) is sufficient to be satisfied by a subset of
$H^{2}(X, \mathbb{Z})$ whose image in $H^{2}(X, \mathbb{Z}_{2})$
is a basis.

Before going on to complex manifolds of dimension~5 we recall that
the Chern character in $K$-theory has a formula:
\begin{equation*}
\ch(\xi) = \hbox{rank} \ \xi + \sum\limits_{k =
1}^{\infty}P_{k}(c_{1}(\xi), \dots, c_{k}(\xi)), \quad \xi\in
K(X),
\end{equation*}
where $P_{k}$ are some universal weighted polynomials with
rational coefficients (see for instance \cite{12}).

In general, for a $r$-tuple of cohomology classes $(c_{1}, \dots,
c_{r}), c_{i}\in H^{2i}(X, \mathbb{Z}), 1\leq i < r$, one defines
\begin{equation*}
\ch(c_{1}, \dots, c_{r}) = r + \sum\limits_{k =
1}^{\infty}P_{k}(c_{1}, \dots, c_{k}),
\end{equation*}
by putting $c_{r + 1} = 0, c_{r + 2} = 0$ and so on. Of course,
when $(c_{1}, \dots, c_{r})$ are the Chern classes of a vector
bundle $E$ of rank $r$, the above expression represents the Chern
character of $E$.

\begin{propo}$\left.\right.$\vspace{.5pc}

\noindent Let $X$ be a connected{\rm ,} compact complex manifold
of complex dimension~$5$. One assumes that the groups $H^{6}(X,
\mathbb{Z})$ and $H^{7}(X, \mathbb{Z})$ have no $2$-torsion{\rm ,}
that $H^{8}(X, \mathbb{Z})$ has no $6$-torsion and that $H^{9}(X,
\mathbb{Z})$ has no $12$-torsion.

Then $(c_{1}, c_{2}, c_{3}, c_{4}, c_{5})$ are the Chern classes
of a {\rm (}unique{\rm )} vector bundle of rank~$5$ on $X$ if and
only if
\begin{equation}
c_{3}\equiv c_{1}c_{2} + \Sq^{2}c_{2} \quad \hbox{in}\ \ H^{6}(X,
\mathbb{Z}_{2}).
\end{equation}
and
\begin{equation}
\int_{X}\ch(c_{1}, \dots, c_{5})\e^{\xi}td(X)\in \mathbb{Z},
\end{equation}
for every $\xi \in H^{2}(X, \mathbb{Z})$.
\end{propo}

\begin{proof}
The cohomology of the manifold $X$ satisfies the assumptions in
Lemma~2.1. Moreover, the coboundaries $d_{2k + 1}^{2p, *}$ of
Atiyah--Hirzebruch spectral sequence vanish for $k > 1$ and $p
\geq 2$. Indeed, $d_{3}^{4}$ vanishes because $2d_{3}^{4} = 0$ and
on the other hand its target, $H^{7}(X, \mathbb{Z})$, has no
2-torsion. Similarly $d_{3}^{6} = 0$ and also $d_{3}^{9} = 0$ by
dimension reasons. Consequently $d_{5}^{4}\hbox{\rm :}\ E_{5}^{4,
0} = H^{4}/\Im d_{3}\rightarrow E_{5}^{9, -4} = H^{9}$. But the
group $H^{9} (X, \mathbb{Z})$ has no 12-torsion and on the other
hand, $12 - d_{5}^{4, *} = 0$ by \cite{7}. Hence $d_{5}^{4} = 0$
as desired.

Let $E$ be a vector bundle on $X$ of rank~5 and let $(c_{1},
\dots, c_{5})$ denote its Chern classes. The congruence~(3.6)
follows from Wu's formula as before. By applying the index theorem
to $E\otimes L$, where $L$ is a line bundle with $c_{1}(L) = \xi$
one gets~(3.7).

Conversely, assume that the cohomology classes $(c_{1}, \dots,
c_{5})$ satisfy~(3.6) and (3.7). In view of Lemma~2.1 we have to
prove that $\dot{c}_{3} = \varphi_{2}(c_{1}, c_{2}), \dot{c}_{4} =
\varphi_{3}(c_{1}, c_{2}, c_{3})$ and $\dot{c}_{5} =
\varphi_{4}(c_{1}, c_{2}, c_{3}, c_{4})$.

We already know the form of the function $\varphi_{2}$ derived
from Wu's formula: $\varphi_{2}(c_{1}, c_{2}) = c_{1}c_{2} +
\Sq^{2}(c_{2})~\,(\bmod\ 2)$. Consequently, there exists a vector
bundle $E$ on $X$ with $c_{i}(E) = c_{i}, i = 1, 2, 3$. Then by
definition, $\varphi_{3}(c_{1}, c_{2}, c_{3})$ is the image of
$c_{4}(E)$ in $H^{8}(X, \mathbb{Z}_{6})$, and thus we are led to
prove that $\overline{\dot{c}_{4}(E)} = \dot{c}_{4}$.

The hypotheses on $H^{*}(X, \mathbb{Z})$ imply
\begin{equation*}
H^{8}(X, \mathbb{Z}_{6})\cong \Hom_{\mathbb{Z}}(H^{2}(X,
\mathbb{Z}), \mathbb{Z}_{6}).
\end{equation*}
Therefore it suffices to prove that $\xi c_{4} - \xi c_{4}(E)$ is
a multiple of 6 for every $\xi \in H^{2}(X, \mathbb{Z})$. The
assumption~(3.7) and the index theorem yield
\begin{equation}
\int_{X}\ch (c_{1}, \dots, c_{5})(1 - \e^{\xi})td (X) \in
\mathbb{Z}
\end{equation}
and
\begin{equation}
\int_{X}\ch (E)(1 - \e^{\xi})td (X) \in \mathbb{Z}
\end{equation}
respectively. \pagebreak

It is worth mentioning that $c_{5}$ and $c_{5}(E)$ do not matter
in the integrality conditions~(3.8) and (3.9). In fact, these
conditions can be reduced modulo $\mathbb{Z}$ to expressions such
as
\begin{align*}
&\frac{1}{6}\xi \cdot c_{4} + P(\xi, c_{1}, c_{2}, c_{3},
c(X)),\\[.4pc]
&\frac{1}{6} \xi \cdot c_{4}(E) + P(\xi, c_{1}, c_{2}, c_{3}
c(X)),
\end{align*}
where $P$ is a polynomial with rational coefficients.

In conclusion, $\xi (c_{4} - c_{4}(E))$ is divisible by~6 for any
$\xi \in H^{2}(X, \mathbb{Z})$, whence $\dot{c}_{4} =
\varphi_{3}(c_{1}, c_{2}, c_{3})$. Therefore, there exists a
vector bundle $F$ on $X$, with $c_{i}(F) = c_{i}, i = 1, 2, 3, 4$.
By simply writing the index theorem for $F$,
\begin{equation*}
\int_{X}\ch (F)td (X)\in \mathbb{Z}
\end{equation*}
and comparing this relation to the assumption~(3.7) one gets that
$c_{4} - c_{4}(F)$ is a multiple of~4!.

This finishes the proof of Proposition~3.2.\hfill $\Box$
\end{proof}

As in the preceding proof, it is sufficient to ask in the
statement of Proposition~3.2 if the class $\xi$ belongs to a
system of generators of $H^{2}(X, \mathbb{Z})$ and $\xi = 0$.

We do not know whether condition~(3.6) is equivalent, like in the
case of $\dim X = 4$, to a congruence involving only the cap
product operation. Also, the stable classification of complex
vector bundles on 6-dimensional manifolds (more exactly a concrete
form of the function $\varphi_{3}$) is unknown to us.

\section{The stable-rank classification on compact complex
manifolds with torsion-free cohomology}

\setcounter{theore}{0}

\setcounter{equation}{0}

In this section we prove that Atiyah--Singer integrality
conditions~(3.3) are the only restrictions imposed to a $n$-tuple
of cohomology classes $(c_{1}, \dots, c_{n})$ in order to be the
Chern classes of a vector bundle on a compact complex manifold
$X$, whose cohomology is torsion free, $n = \dim_{C}X$.

\begin{propo}$\left.\right.$\vspace{.5pc}

\noindent Let $X$ be a connected{\rm ,} compact complex manifold
of dimension $n${\rm ,} with torsion-free cohomology. Then the
cohomology classes $(c_{1}, \dots, c_{n})$ are the Chern classes
of a {\rm (}unique{\rm )} vector bundle of rank $n$ on $X${\rm ,}
if and only if \setcounter{equation}{0}
\begin{equation}
\int_{X}\ch (c_{1}, \dots, c_{n})\ch (\xi)td (X)\in
\mathbb{Z},\quad \xi \in K(X).
\end{equation}
\end{propo}

\begin{proof}
Conditions~(4.1) are necessary by the Atiyah--Singer index
theorem. Conversely, let $(c_{1}, \dots, c_{n})$ be an $n$-tuple
of cohomology classes which satisfies~(4.1). In order to construct
a vector bundle $E$ with $c_{i}(E) = c_{i}, 1\leq i\leq n$, it is
sufficient to prove that
\begin{equation*}
\dot{c}_{3} = \varphi_{2}(c_{1}, c_{2}), \dots,\quad \dot{c}_{n} =
\varphi_{n - 1}(c_{1}, \dots, c_{n - 1})
\end{equation*}
and apply them to Lemma~2.1.

However, in this case Lemma~2.1 can be avoided and a direct
argument is possible. So we prove inductively that there are
vector bundles $E^{p}$ on $X$, with the properties: $c_{i}(E^{p})
= c_{i}, 1\leq i\leq p$ and $1\leq p\leq n$.

The case $p = 1$ is obvious. Assume that there exists a vector
bundle $E^{p - 1}$ with $c_{i}(E^{p - 1}) = c_{i}$ for $1\leq
i\leq p - 1$. We have to prove that $c_{p} - c_{p}(E^{p - 1})$ is
a multiple of $(p - 1)!$. Indeed, in that case the result of
Atiyah and Hirzebruch mentioned in the introduction provides a
vector bundle $F$ with $c_{i}(F) = 0$ for $1\leq i\leq p - 1$ and
$c_{p}(F) = c_{p} - c_{p}(E^{p - 1})$. Hence the vector bundle
$E^{p} = E^{p - 1}\oplus F$ has the first Chern classes
$c_{i}(E^{p}) = c^{i}, 1\leq i\leq p$.

Since
\begin{equation*}
H^{2p}(X, \mathbb{Z}_{(p - 1)!})\cong \Hom_{\mathbb{Z}} (H^{2n
-2p}(X, \mathbb{Z}), \mathbb{Z}_{(p - 1)!}),
\end{equation*}
it is sufficient to verify that $\eta(c_{p} - c_{p}(E^{p - 1}))$
is a multiple of $(p - 1)!$ for any $\eta \in H^{2n - 2p}(X,
\mathbb{Z})$. For that aim, consider a class $\xi \in
\tilde{K}(X)$ with $c_{i}(\xi) = 0$ for $0\leq i < n - p$ and
$c_{n - p}(\xi) = (-1)^{n - p - 1}(n - p - 1)!\eta$, so that
\begin{equation*}
\ch (\xi) = \eta + \hbox{higher order terms}.
\end{equation*}
By hypothesis
\begin{equation*}
\int_{X}\ch(c_{1}, \dots, c_{n})\ch(\xi)td(X)\in \mathbb{Z},
\end{equation*}
and by the index theorem,
\begin{equation*}
\int_{X}\ch(E^{p - 1}\ch (\xi)td (X)\in \mathbb{Z}.
\end{equation*}
In both integrals the classes $c_{p + 1}, \dots, c_{n}$ and $c_{p
+ 1}(E^{p - 1}), \dots, c_{n}(E^{p - 1})$ do not matter.
Consequently, they reduce to
\begin{equation*}
\frac{1}{(p - 1)!}\eta c_{p} + \mathbb{Q}(c_{1}, \dots, c_{p - 1},
c(\xi), c(X))\in \mathbb{Z}
\end{equation*}
and
\begin{equation*}
\frac{1}{(p - 1)!}\eta c_{p}(E^{p - 1}) + \mathbb{Q}(c_{1}, \dots,
c_{p - 1}, c(\xi), c(X))\in \mathbb{Z}
\end{equation*}
respectively, where $\mathbb{Q}$ is a polynomial with rational
coefficients.

In conclusion, $(p - 1)!$ divides $\eta(c_{p} - c_{p}(E^{p - 1}))$
and the proof of Proposition~4.1 is complete.\hfill $\Box$
\end{proof}

Before continuing with an application of the proposition we recall
some needed terminology. A~nonsingular projective (complex)
variety $X$ of dimension $n$ is said to possess an algebraic cell
decomposition
\begin{equation*}
X = X^{n}\supset X^{n - 1}\supset \cdots \supset X^{0}\supset
X^{-1} = \varnothing,
\end{equation*}
if $X^{i}$ are closed algebraic subsets of $X$ and each difference
$X^{i}\backslash X^{i - 1}$ is a disjoint union of locally closed
submanifolds $U^{ij}$ which are isomorphic with the affine space
$C^{i}$. The sets $U_{ij}$ are the cells of this decomposition and
their closures are called the closed cells.

If the projective manifold $X$ admits an algebraic cell
decomposition, then $H^{\rm odd}(X, \mathbb{Z}) = 0, H^{\rm
even}(X, \mathbb{Z})$ is torsion free and moreover $H^{2p}(X,
\mathbb{Z})$ is freely generated by the fundamental classes of all
closed cells on codimension $p$ (see \cite{9}). Concerning the
$K$-theory, it is known that $K_{\rm alg}(X)$ is a free group with
base $[{\mathcal O}_{Y}]$, where ${\mathcal O}_{Y}$ are the
structure sheaves of the closed cells (see \cite{11}). Moreover,
the natural map
\begin{equation}
\varepsilon\hbox{\rm :}\ K_{\rm alg}(X)\rightarrow K(X)
\end{equation}
is in that case isomorphic.

The surjectivity of $\varepsilon$ follows easily from the
following statement: if $\xi \in K(X)$ has the property
$c_{i}(\xi) = 0$ for $1 \leq i \leq p - 1$ then there exists an
algebraic class $\xi'\in K_{\rm alg}(X)$ with $c_{i}(\xi) =
c_{i}(\varepsilon(\xi'))$ for $1 \leq i \leq p$.

If $\xi \in K(X)$ is as before, then $c_{p}(\xi) = (-1)^{p - 1}(p
- 1)!\eta$, with $\eta \in H^{2p}(X, \mathbb{Z})$. Since $\eta -
\sum_{i}n_{i}[Y_{i}]$, where $Y_{i}$ are the closed cells of
dimension $p$, we define $\xi' = \sum_{i}n_{i}[{\mathcal
O}_{Y_{i}}]$. But for an irreducible subvariety $Y$ of codimension
$p$ one has $c_{i}({\mathcal O}_{Y}) = 0$ for $0 < i < p$ and
$c_{p}({\mathcal O}_{Y}) = (-1)^{p - 1}(p - 1)![Y]$ (see p.~297 of
\cite{9} for an argument valid when $Y$ is nonsingular) (in
general, by simply removing the singular locus of $Y$). This
concludes the surjectivity. The injectivity of $\varepsilon$
follows\break similarly.

\begin{propo}$\left.\right.$\vspace{.5pc}

\noindent Let $X$ be a projective{\rm ,} nonsingular variety of
dimension $n${\rm ,} which has an algebraic cell decomposition.
The cohomology classes $(c_{1}, \dots, c_{n})$ are the Chern
classes of a {\rm (}unique{\rm )} vector bundle of rank $n$ on $X$
if and only if
\begin{equation}
\int_{X}\ch (c_{1}, \dots, c_{n})\ch ({\mathcal O}_{Y})td (X)\in
\mathbb{Z},
\end{equation}
for every closed cell $Y$ of $X${\rm ,} of dimension greater than
or equal to~$3$.
\end{propo}

\begin{proof}
The necessity follows from the index theorem or from the
Riemann--Roch--Hirzebruch formula.

If we drop the condition $\dim Y \geq 3$, the statement is
equivalent to Proposition~4.1. So the only thing to be proved is
that the closed cells $Y$ of dimension less than~3 do not matter
in the integrality formula~(4.3). This is obvious when we compare
with the proof of Proposition~4.1, because any vector bundle $E$
on $X$ with $c_{i}(E) = c_{i}, 1 \leq i < p$ satisfies $c_{p} =
c_{p}(E)$ in $H^{2p}(X, \mathbb{Z}_{(p - 1)!})$ for $p \leq
2.$\hfill $\Box$
\end{proof}

Otherwise, one repeats the proof of Proposition~4.1 by evaluating
$c_{p} - c_{p}(E)$ on $\eta = [Y]$ and by taking $\xi =
\varepsilon [{\mathcal O}_{Y}]$, correspondingly.

\begin{rema}
{\rm Let $Y$ be a connected, nonsingular, closed algebraic subset
of dimension $d$ of a projective nonsingular variety $X$ of
dimension $n$.

The Riemann--Roch--Grothendieck theorem applied to the inclusion
$i\hbox{\rm :}\ Y\hookrightarrow X$\break yields
\begin{equation*}
i_{*}(td(Y)) = \ch ({\mathcal O}_{Y})td(X).
\end{equation*}
This identity and the projection formula prove that for any
cohomology classes $(c_{1}, \dots, c_{n})$ one has}
\begin{align*}
\int_{X}\ch (c_{1}, \dots, c_{n})\ch({\mathcal O}_{Y})td(X) &=
\int_{Y}\ch(i^{*}c_{1}, \dots, i^{*}c_{n})td(Y)\\[.4pc]
&= \int_{Y}\ch(i^{*}c_{1}, \dots, i^{*}c_{d})td(Y)\\[.4pc]
&\quad\,+ (n - d)\chi ({\mathcal O}_{Y})
\end{align*}
because $i^{*}(c_{d+1}) = \cdots = i^{*}(c_{n}) = 0$.
\end{rema}

This remark shows that condition~(4.3) in Proposition~4.2 can be
replaced by the Riemann--Roch formula, written on every closed
cell, provided that all closed cells are non-singular. A~familiar
example of that kind is $X = \mathbb{P}^{n}(\mathbb{C})$, the
complex projective space of dimension $n$. Actually Le Potier
\cite{16} had a similar point of view on
$\mathbb{P}^{n}(\mathbb{C})$, in what concerns the integrability
conditions satisfied by the Chern classes of a vector bundle.

Quite specifically, in that case $(c_{1}, \dots, c_{n})$ are the
Chern classes of a vector bundle on $\mathbb{P}^{n}(\mathbb{C})$,
if and only if
\begin{equation*}
\int_{\mathbb{P}^{3}}\ch(c_{1}, c_{2}, c_{3})td
(\mathbb{P}^{3})\in \mathbb{Z}, \dots, \int_{\mathbb{P}^{n}}\ch
(c_{1}, \dots, c_{n})td (\mathbb{P}^{n})\in \mathbb{Z}.
\end{equation*}
Recall that $td(\mathbb{P}^{n}) = [h/(1 - \hbox{e}^{-h})]^{n +
1}$, where $h$ is the class of a hyperplane.

An equivalent set of integrality conditions, such as
\begin{equation*}
\int_{\mathbb{P}^{n}}\ch (c_{1}, \dots, c_{n})\hbox{e}^{kh}[h/(1 -
\e^{-h})]^{n + 1}\in \mathbb{Z}, \quad 0\leq k\leq n - 3,
\end{equation*}
is provided by the next result.

\begin{propo}$\left.\right.$\vspace{.5pc}

\noindent Let $X$ be a complex compact manifold of dimension $n$.
Assume that $H^{*}(X, \mathbb{Z})$ is torsion-free and that
$H^{\rm even}(X, \mathbb{Z})$ is generated multiplicatively by
$H^{2}(X, \mathbb{Z})$. Then $(c_{1}, \dots, c_{n})$ are the Chern
classes of a {\rm (}unique{\rm )} vector bundle of rank $n$ on $X$
if and only if
\begin{equation}
\int_{X}\ch(c_{1}, \dots, c_{n})\e^{\xi}td(X)\in \mathbb{Z},\quad
\xi \in H^{2}(X, \mathbb{Z}).
\end{equation}
Moreover{\rm ,} the class $\xi$ can be chosen of the form $\xi =
\xi_{j_{1}} + \cdots + \xi_{j_{k}}${\rm ,} where $0\leq k\leq n -
3, 1\leq j_{1}\leq \cdots \leq j_{k}\leq m$ and $\xi_{1}, \dots,
\xi_{m}$ is a basis of $H^{2}(X, \mathbb{Z})$ {\rm (}we put $\xi =
0$ for $k = 0${\rm )}.
\end{propo}

\begin{proof}
Conditions~(4.4) are necessary by the Atiyah--Singer index
theorem.

Conversely, let $(c_{1}, \dots, c_{n})$ be the cohomology classes
which satisfy~(4.4). For proving the existence of a vector bundle
on $X$ with the prescribed Chern classes $(c_{1}, \dots, c_{n})$
one may use Lemma~2.1 or one may argue as in the proof of
Proposition~4.1. As a matter of fact, we have to prove that for
any $p \geq 3$ and any vector bundle $E$ with $c_{i}(E) = c_{i}$
for $1\leq i < p$, the classes $c_{p}$ and $c_{p}(E)$ have the
same image in $H^{2p}(X, \mathbb{Z}_{(p - 1)!})$ or equivalently,
$\eta(c_{p} - c_{p}(E))$ is divisible by $(p - 1)!$ for a set of
generators $\eta$ of $H^{2n -2p}(X, \mathbb{Z})$.

We may assume that $\eta = \zeta_{1}, \dots, \zeta_{n - p}$, where
$\zeta_{i} = \xi_{j_{i}}$ for $1\leq i\leq n - p$ and certain
indices $j_{1}, \dots, j_{n - p}$. Let us denote
\begin{align*}
K(\zeta_{1}, \dots, \zeta_{n - p}) = 1 -
\sum\limits_{i}\e^{\zeta_{i}} + \sum\limits_{j\neq k}\e^{\zeta_{j}
+ \zeta_{k}} + \cdots + (-1)^{n - p}\e^{\zeta_{i} + \cdots +
\zeta_{n - p}},
\end{align*}
so that assumption~(4.4) becomes
\begin{equation}
\int_{X}\ch(c_{1}, \dots, c_{n})K(\zeta_{1}, \dots, \zeta_{n -
p})td(X)\in \mathbb{Z}.
\end{equation}

On the other hand, the index theorem applied to bundles of the
form $E\otimes L, \hbox{rk}\,L = 1$, yields
\begin{equation}
\int_{X}\ch(E)K(\zeta_{1}, \dots, \zeta_{n - p})td(X)\in
\mathbb{Z}.
\end{equation}
Next it is sufficient to remark that
\begin{align*}
K(\zeta_{1}, \dots, \zeta_{n - p}) &= (1 - \e^{\zeta_{1}})\cdots
(1 - \e^{\zeta_{n - p}})\\[.4pc]
&= (-1)^{n - p}\zeta_{1}\dots \zeta_{n - p} + \hbox{higher order
terms}\\[.4pc]
&= (-1)^{n - p}\eta + \hbox{higher order terms}.
\end{align*}
This shows that the classes $c_{p + 1}, \dots, c_{n}$ and $c_{p +
1}(E), \dots, c_{n}(E)$ do not matter in relations~(4.5) and
(4.6), and moreover by subtracting~(4.6) from (4.5) one finds that
$(p - 1)!$ divides $\eta(c_{p} - c_{p}(E))$.

This completes the proof of Proposition~4.3.\hfill $\Box$
\end{proof}

\section{The stable-rank classification on flag manifolds}

\setcounter{theore}{0}

\setcounter{equation}{0}

A notable class of algebraic manifolds which fit into the
conditions of Proposition~4.2 are the flag manifolds, or
generalized flag manifolds, following the terminology adopted by
some authors.

Let $G$ denote a semi-simple, connected and simply-connected
complex algebraic linear group. One denotes by $T$ a maximal torus
of $G$, by $N(T)$ the normalizer of $T$ in $G$ and by $W = N(T)/T$
the corresponding Weyl group. Let $B$ be a Borel subgroup of $G$
which contains $T$, $\Phi$ the root system with respect to $T$ and
$\Phi^{+}\subset \Phi$ the positive roots corresponding to the
choice of $B$. The system of simple (positive) roots is denoted by
$\Delta \subset \Phi^{+}$.

The parabolic subgroup $P_{I} = BW_{I}B$ associated to a subset $I
\subset \Delta$ is obtained from the subgroup $W_{I}\subset W$
generated by the reflections corresponding to $I$. Any parabolic
subgroup of $G$ can be obtained by such a standard construction
(see for instance \cite{5}).

Let $P = P_{I}$ be a parabolic subgroup of $B$. We simply denote
$W_{P}$ instead of $W_{I}$. Let $w\in W$ and let $n(w)$ denote a
representative of $w$ in $N(T)$. One easily remarks that the
double $(B, P)$ coset $Bn(w)P$ depends only on $wW_{P}$ and hence
it may be denoted by $BwP$. This is the so-called open {\it Bruhat
cell} associated to $wW_{P}$. The images of the Bruhat cells
through the natural projection $\pi\hbox{\rm :}\ G \rightarrow
G/P$ give an algebraic cell decomposition of the {\it flag
manifold} $G/P$. The corresponding closed cells, called the {\it
Schubert cells} of $G/P$, are denoted by $X_{P}(w)$, and they are
labelled as $w \in W/W_{P}$. In fact, $B$ acts (from the left) on
$G/P$ and the Schubert cells are exactly the closures of the
corresponding orbits.


If $I$ is empty, then $P = B$ and the corresponding homogeneous
space $G/B$ is called the full flag manifold of $G$. In that case
the Schubert cells $X_{w}$ are labelled by the Weyl group $W$.

The flag manifolds are projective rational homogeneous varieties,
and conversely any such abstract variety is isomorphic to a flag
manifold (see \cite{5}).

Recall now the Atiyah--Hirzebruch morphism for the full flag
manifold $G/B$ \cite{1}. Let $\hat{T}$ be the group of characters
of $T$ and let $R(T) = Z[\hat{T}]$ be the representation ring of
$T$. Any character $\chi \in \hat{T}$ extends trivially over the
unipotent radical of $B$ and it gives a complex representation
$\chi\hbox{\rm :}\ B\rightarrow C^{*}$. Thus one obtains a line
bundle $L(\chi) = G\times_B C^{*}$ on $G/B$. The association
$\chi\mapsto L(\chi)$ is a morphism of groups $\hat{T}\rightarrow
\hbox{Pic}(G/B)$ which extends to a morphism of rings
$R(T)\rightarrow K(X)$, known as the Atiyah--Hirzebruch morphism.

Since the group $G$ was supposed to be simply connected, the map
$\hat{T} \rightarrow\hbox{Pic}(G/B)$ is actually an isomorphism.
The first Chern class gives an isomorphism $\hbox{Pic}(G/B) \cong
H^{2}(G/B, \mathbb{Z})$. For simplicity, we identify the integer
weights of $G$ with the characters of $T$ and with elements of
$H^{2}(G/B, \mathbb{Z})$.

Consider a flag manifold $X = G/P$ of dimension $n$, as above.
Proposition~4.2 tells us that a system of cohomology classes
$(c_{1}, \dots, c_{n})$ of $X$ is the total Chern class of a
vector bundle on $X$ if and only if
\begin{equation*}
\int_{X} \ch (c_{1},\dots, c_{n}) \ch ({\mathcal O}_{X_p(w)})td(X)
\in \mathbb{Z}, \quad w \in W/W_{P}.
\end{equation*}
Our next aim is to lift these conditions on the full flag manifold
$G/B$. Let $\pi\hbox{\rm :}\ G/B \rightarrow G/P$ be the canonical
projection. This is a fibration with fibre the rational manifold
$P/B$. The induced morphism $\pi^{*}\hbox{\rm :}\ H^{*} (G/P,
\mathbb{Z})\rightarrow H^{*} (G/B, \mathbb{Z})$ is one-to-one and
identifies the cohomology of $G/P$ with the ring of
$W_{P}$-invariants in $H^{*} (G/B, \mathbb{Z})$ \cite{4}.

Let $w \in W$. One has
\begin{equation}
\pi_{*} ({\mathcal O}_{X_{w}}) = {\mathcal O}_{X_{P}(w)};\quad
R^{q}\pi_{*} ({\mathcal O}_{X_{w}}) = 0,\quad q \geq 1.
\end{equation}
This is a subtle result essentially proved by Ramanan and
Ramanathan \cite{19} (see also Proposition~2 of \S5 in \cite{8}).
Consequently $\pi_{*}[{\mathcal O}_{X_{w}}] = ({\mathcal
O}_{X_{P}(w)})$ at the level of algebraic $K$-theories.

Riemann--Roch--Grothendick theorem and the projection formula
yield
\begin{align}
&\pi_{*}[\ch (c_{1},\dots, c_{n}) \ch ({\mathcal O}_{X_{w}}) td
(G/B)]\nonumber\\[.2pc]
&\quad\, = \ch (c_{1},\dots, c_{n}) \pi_{*}[\ch({\mathcal
O}_{X_{w}}) td
(G/B)]\nonumber\\[.2pc]
&\quad\, = \ch (c_{1},\dots,c_{n}) \ch ({\mathcal
O}_{X_{P}(w)})td(G/P).
\end{align}
Recall that $td(G/B) = \e^{-\rho}$, where $\rho$ is the semisum of
all positive roots of $G$ (see \cite{6}).

By evaluating formula (5.2) on $G/P$ one obtains
\begin{align*}
&\int_{G/B} \ch (c_{1},\dots, c_{n}) \ch ({\mathcal O}_{X_{w}})
\e^{-\rho}\\[.2pc]
&\quad\, = \int_{G/P} \ch (c_{1},\dots, c_{n}) \ch ({\mathcal
O}_{X_{P}(w)}) td(G/P).
\end{align*}
In particular, this shows that the left-hand term does not depend
on the representative $w$ of the class $wW_{P} \in W/W_{P}$ and
explains the meaning of (5.3).

On the other hand, it is known that the Atiyah--Hirzebruch map is
onto (see \cite{18}). Therefore any class $\ch (\xi)$ with $\xi
\in K(G/B)$ is an entire combination of elements of the form
$\e^\tau$ with $\tau \in \hat{T} = H^{2}(G/B, \mathbb{Z})$. Hence
we have proved the following.

\begin{theor}[\!]
Let $X = G/P$ be a flag manifold of dimension $n$. The $n$-tuple
of cohomology classes $(c_{1},\dots,c_{n})$ of $X$ represents the
Chern classes of a {\rm (}unique{\rm )} vector bundle of rank $n$
on $X$ if and only if
\begin{equation}
\int_{G/B} \ch(c_{1},\dots,c_{n}) \ch ({\mathcal O}_{X_{w}})
\e^{-\rho} \in \mathbb{Z} \quad for \ \ w \in W/W_{P},
\end{equation}
or equivalently{\rm ,} if and only if
\begin{equation}
\int_{G/B} \ch(c_{1},\dots,c_{n}) \e^{\chi} \in \mathbb{Z} \quad
for \ \ \chi \in \hat{T}.
\end{equation}
\end{theor}

\begin{remarr}$\left.\right.$

\begin{enumerate}
\renewcommand\labelenumi{(\arabic{enumi})}
\leftskip .15pc
\item According to Proposition~4.2, it is sufficient to consider in
relation (5.3) only cosets $wW_{P}$ with the property that $\dim
X_{P}(w) \geq 3$.

\parindent 1pc Let $Y$ be a Schubert cell in $X$. Then $\pi^{-1}(Y)$ is an
irreducible $B$-invariant subset of $G/B$, so it coincides with a
Schubert cell of $G/B$. Thus for any coset $wW_{P}$ one may find a
representative $w \in W$, such that $\pi^{-1}(X_{P}(w)) = X_{w}$
($w$ is the element of maximal length in $wW_{P}$). Thus $w \in W$
is called $P$-saturated. In that case the vanishing conditions
(5.1) easily follow from the fact that $\pi/X_{w}\hbox{\rm :}\
X_{w} \rightarrow X_{P}(w)$ is a fibration of fiber $P/B$, and
$H^{q}(P/B, {\mathcal O}\;) = 0$ for $q\geq 1$. Thus one can avoid
ref.~\cite{19} and the corresponding integrality conditions (which
replace (5.3)) are
\begin{equation*}
\hskip -1.25pc \int_{G/B} \ch (c_{1},\dots, c_{n}) \ch ({\mathcal
O}_{X_{w}}) \e^{-\rho} \in \mathbb{Z},
\end{equation*}
for any $P$-saturated word $w \in W$ of length $l(w) \geq \dim
(P/B) + 3$.

\item As $\int_{G/B} \ch(c_{1},\dots, c_{n}) \e^{\tau} =
\int_{G/B} \ch(c_{1},\dots, c_{n}, 0,\dots, 0) \e^{\tau} +$ and
integer, the integrality conditions (5.3) shows that the
classification of stable-rank vector bundles on $G/P$ can be
deduced from the same classification on $G/B$.

\parindent 1pc Namely $\pi^{*}\hbox{\rm :}\ \Vect_{\rm top}^{s,{\rm rk}} (G/P)
\rightarrow \Vect_{\rm top}^{s,{\rm rk}} (G/B)$ is a one-to-one
map and a vector bundle $E$ on $G/B$ (of stable-rank) is in the
range of $\pi^{*}$ if and only if the Chern classes $c_{i}(E)$ are
$W_{p}$-invariant when $i \geq 1$. Of course, the two stable ranks
on $G/P$ and $G/B$ do not coincide and the precise meaning of
$\pi^{*}(E)$ is up to a trivial direct summand.

\parindent 1pc By keeping the notation introduced in the first part of this
section we denote by $w_{0}$ in the sequel the element of maximal
length in $W$. Also $w_{1},\dots, w_{r}$ are the fundamental
weights of $G$ with respect to the Borel subgroup $B$.
\end{enumerate}
\end{remarr}

\begin{theor}[\!]
Let $X = G/P$ be a flag manifold and assume that $H^{2}(X,
\mathbb{Z})$ generates $H^{\rm even} (X, \mathbb{Z}), n = \dim X$.
Then $\Vect_{\rm top}^{n}(X)$ is in bijection{\rm ,} via Chern
classes{\rm ,} with the $n$-tuples $(c_{1},\dots, c_{n})$
satisfying
\begin{equation}
\int_{G/B} \ch (c_{1},\dots, c_{n}) \e^{\chi-\rho} \in \mathbb{Z},
\end{equation}
for every character $\chi \in \hat{T}$ of the form $\chi =
\omega_{i_{1}} + \cdots + \omega_{i_{k}}${\rm ,} where $0 \leq k
\leq n - 3, 1\leq i_{1} \leq \dots \leq i_{k} \leq r$ and
$\omega_{i_{1}}, \dots, \omega_{i_{k}}$ are the fundamental
weights corresponding to simple roots $\alpha$  for which
$w_{0}\sigma_{\alpha}$ is $P$-saturated.
\end{theor}

\begin{proof}
One applies Proposition~4.3. Let $\alpha$ be a root with
$w_{0}\sigma_{\alpha}$ $P$-saturated ($\sigma_{\alpha}$ is the
associated reflection). Then $X_{P} (w_{0}\sigma_{\alpha})$ is a
hypersurface in $X$, whose corresponding line bundle $L_{P}
(\alpha) = {\mathcal O}_{X}(X_{P}(w_{0}\sigma_{\alpha}))$ has the
first Chern class $[X_{P} (w_{0}\sigma_{\alpha})]$.

By varying $\alpha$ one gets a basis for $\hbox{Pic}(X) \cong
H^{2} (X, \mathbb{Z})$. On the other hand,
\begin{equation*}
\pi^{*} L_{P} (\alpha) = {\mathcal O}_{G/B}
(X_{w_{0}\sigma_{\alpha}}) = L(\omega_{\alpha}),
\end{equation*}
where $\omega_{\alpha}$ is the fundamental weight associated to
$\alpha$. Since for any line bundle $L$ on $X$ one has
\begin{equation*}
\int_{X} \ch (c_{1},\dots, c_{n}) \ch (L) td (X) = \int_{G/B} \ch
(c_{1},\dots, c_{n}) \ch (\pi^{*} L) \e^{-\rho},
\end{equation*}
the proof of Theorem~5.2 is over.\hfill $\Box$
\end{proof}

Let $S(\hat{T})$ denote the symmetric algebra of $\hat{T}$, i.e.
the polynomials with integer coefficients acting on the Lie
algebra of $T$. The isomorphism
\begin{equation*}
\hat{T} \cong H^{2}(G/B, \mathbb{Z}), \quad \chi \mapsto
c_{1}(L(\chi)),
\end{equation*}
extends to a morphism of rings
\begin{equation*}
c\hbox{\rm :}\ S(\hat{T}) \rightarrow H^{*}(G/B, \mathbb{Z}).
\end{equation*}
Demazure \cite{8} proved that $c$ is onto precisely when $G$ is a
product of groups of type SL and Sp. In particular in that case
$H^{2}(G/B, \mathbb{Z})$ generates $H^{\rm even}(G/B,
\mathbb{Z})$, and so this fits into the conditions of Theorem~5.2.
However, for a parabolic subgroup $P \subset G, H^{2}(G/P,
\mathbb{Z})$ may not generate the even cohomology of $G/P$. In
spite of that, the next result still holds.

\begin{theor}[\!]
Let $X = G/P$ be an $n$-dimensional flag manifold with $G$ a
product of simple Lie groups of types {\rm SL} and {\rm Sp}.

Then $\Vect_{\rm top}^{n}(X)$ is in bijection{\rm ,} via Chern
classes{\rm ,} with the set of $n$-tuples $(c_{1},\dots, c_{n}),$
$c_{i} \in H^{2i}(G/P, \mathbb{Z}), 1 \leq i \leq n${\rm ,}
satisfying
\begin{equation*}
\int_{G/B} \ch (c_{1},\dots, c_{n}) \e^{x-\rho} \in \mathbb{Z},
\end{equation*}
where $\chi = \omega_{i_{1}} + \cdots + \omega_{i_{k}}, 0\leq k
\leq \dim(G/B)-3, 1\leq i_{1} \leq \dots \leq i_{k} \leq r$.
\end{theor}

\begin{proof}
If $P = B$, then the statement reduces to Theorem~5.2. Otherwise
one uses the equality $\ch (c_{1}, \dots, c_{n}) = \ch
(c_{1},\dots, c_{n}, 0,\dots, 0) + n-\dim(G/B)$ and Remark~2.

We close this section by a few comments on the Lie algebra
interpretation of the integrality conditions. The morphism $c$
induces a ring morphism
\begin{equation*}
c_{\mathbb{Q}}\hbox{\rm :}\ S_{\mathbb{Q}}(\hat{T}) \rightarrow
H^{*} (G/B, \mathbb{Q}),
\end{equation*}
where $S_{\mathbb{Q}}(\hat{T}) =
\mathbb{Q}\otimes_{\mathbb{Z}}S(\hat{T})$. The Weyl group $W$ acts
on $S_{\mathbb{Q}}(\hat{T})$ and a classical result asserts that
$c_{\mathbb{Q}}$ is onto and ker $c_{\mathbb{Q}}$ is the ideal
generated by the homogeneous $W$-invariant polynomials of positive
degree \cite{5}.

For any $w \in W$ one defines the operators
\begin{equation*}
A_{w}\hbox{\rm :}\ S_{\mathbb{Q}}(\hat{T}) \rightarrow
S_{\mathbb{Q}}(\hat{T}), \quad D_{w}\hbox{\rm :}\
S_{\mathbb{Q}}(\hat{T}) \rightarrow \mathbb{Q}
\end{equation*}
by $A_{\alpha}f = (f-\sigma_{\alpha}f)/\alpha, f\in
S_{\mathbb{Q}}(\hat{T})$ for a simple root $\alpha, A_{w} =
A_{\alpha_{1}} \dots A_{\alpha_{l}}, D_{w} f = (A_{w}f)(0)$ for
any irreducible decomposition $w = \sigma_{\alpha_{1}}\dots
\sigma_{\alpha_{l}}$ (see \cite{4}).

The operators $A_{w}$ and $D_{w}$ descend to the rational
cohomology of $G/B$ (we use, for convenience, the same letters)
\begin{equation*}
A_{w}\hbox{\rm :}\ H^{*}(G/B, \mathbb{Q}) \rightarrow H^{*} (G/B,
\mathbb{Q}),\quad D_{w}\hbox{\rm :}\ H^{*} (G/B,
\mathbb{Q})\rightarrow \mathbb{Q}.
\end{equation*}

According to Bermstein--Gelfand--Gelfand \cite{4}, the integrality
condition (5.4) can be expressed as
\begin{equation*}
\int_{G/B} \ch (c_{1},\dots, c_{n}) \e^{\tau} = D_{w_{0}}(\ch
(c_{1},\dots, c_{n}) \e^{\tau}),
\end{equation*}
for any $\tau \in \hat{T}$. Recall that $w_{0}$ stands for the
maximal length element in $W$.

Concerning the integrality conditions (5.3), their translation
into Lie algebra terms faces an additional difficulty. More
exactly, $H^{*}(G/B, \mathbb{Q})$ has two natural bases as
$\mathbb{Q}$-vector space
\begin{equation*}
([X_{w}])_{w \in W} \quad\hbox{and}\quad (\ch ({\mathcal
O}_{X_{w}}))_{w \in W},
\end{equation*}
related by some linear formulae:
\begin{equation*}
\ch ({\mathcal O}_{X_{w}}) = \sum\limits_{w'\in W}
q_{w,w'}[X_{w'}], \quad w \in W,
\end{equation*}
where $q_{w,w'} \in \Q, q_{w,w'}=0$ unless $w'\leq w$ and $q_{w,w}
=1$. Therefore, by \cite{4}, the expression in (5.4) becomes
\begin{align*}
\int_{G/B} \ch(c_1,\dots,c_n) \ch({\mathcal O}_{X_w}) \e^{-\rho} =
\sum\limits_{w'\in W} q_{w,w'} D_{w'} (\ch(c_1,\dots,c_n)
\e^{-\rho}).
\end{align*}

$\left.\right.$\vspace{-1pc}

\hfill $\Box$
\end{proof}

The next natural question is how to compute the matrix
$(q_{w,w'})_{w, w'\in W}$. The results of Demazure \cite{8} on the
desingularization of the Schubert cells of the flag manifolds give
an algorithm for computing this matrix. Unfortunately we were not
able to put this algorithm into a simple, manageable form (see the
subsequent Appendix).

Another solution to this problem would be to find for every $w\in
W$ a concrete polynom in $S_{\mathbb{Q}}(\hat{T})$ whose image
through $c_{\mathbb{Q}}$ is $\ch({\mathcal O}_{X_w})$.

\section*{Appendix: The Chern character of Schubert cells}

\renewcommand\theequation{A\arabic{equation}}

\setcounter{equation}{0}

Let $T\subset B \subset G, \Delta \subset \Phi^+ \subset \Phi$ be
as in \S5 and let $n=\dim(G/B)$. Throughout this appendix,
$\tilde{B}$ stands for the opposite Borel subgroup of $B$ and $w_0
= s_{\beta_1} \dots s_{\beta_n}$ is the element of $W$ of maximal
length, written in a reduced form. One denotes
\begin{equation*}
\alpha_1 = \beta_1, \alpha_2 = s_{\beta_1} (\beta_2),\dots,
\alpha_n = s_{\beta_1} \dots s_{\beta_{n-1}} (\beta_n).
\end{equation*}
Then $\Phi^+ = \{ \alpha_1, \dots,\alpha_n\}$ and $s_{\alpha_1}
\dots s_{\alpha_i} = s_{\beta_i} \dots s_{\beta_1}$ for $1 \leq i
\leq n$.\pagebreak

Let $\psi\hbox{\rm :}\ Z\rightarrow X = G/B$ be the
desingularization associated to the decomposition of $w_0$ into
$s_{\beta_1}\dots s_{\beta_n}$ \cite{8}. The space $Z$ is
constructed as a sequence of $\Pa^1$-fibre bundles $f_i\hbox{\rm
:}\ Z^i \rightarrow Z^{i-1}$, each of them being endowed with a
canonical cross-section $\sigma_i\hbox{\rm :}\ Z^{i-1} \rightarrow
Z^i$,\vspace{-.5pc}
\begin{equation*}
Z = Z^n \begin{array}{c} {_{f_n}}\\[-.4pc]
\longrightarrow\\[-.7pc] \longleftarrow\\[-.9pc] {_{\sigma_n}}\\[.5pc]
\end{array} Z^{n-1} \begin{array}{c} \longrightarrow\\[-.7pc]
\longleftarrow \end{array} \cdots \begin{array}{c} \longrightarrow\\[-.7pc]
\longleftarrow \end{array} Z^1 \begin{array}{c} {_{f_1}}\\[-.4pc]
\longrightarrow\\[-.7pc] \longleftarrow\\[-.9pc] {_{\sigma_1}}\\[.5pc]
\end{array} Z^0 = \hbox{point}.
\end{equation*}

$\left.\right.$\vspace{-1.2pc}

For any $i =1,2,\dots, n$ one denotes $Z_i = f_n^{-1} \dots
f_{i+1}^{-1} (\Im \sigma_i)$ and one remarks that $\Z_i$ are
smooth hypersurfaces in $Z$. Following \cite{8}, each intersection
$Z_K = \bigcap_{i\in K} Z_i, K \subset [1,n]$ is still smooth and
of codimension $|K|$. Let $\xi_i$ denote the fundamental class
$[Z_i] \in H^2(Z,\Z), i = 1,\dots, n$. Consequently $[Z_K] =:
\xi_K = \prod_{i\in K} \xi_i$.

Demazure proves that the cohomology ring $H^* (Z,\Z)$ is
$\Z$-free, with generators $(\xi_K)_{K\subset [1,n]}$ and
relations
\begin{align*}
&\xi_1^2 = 0, \xi_2^2 + \langle \alpha_1^v, \alpha_2 \rangle \xi_1
\xi_2 = 0, \dots,\\[.2pc]
&\xi_n^2 + \langle \alpha_1^v, \alpha_n \rangle \xi_1 \xi_n +
\dots + \langle \alpha_{n-1}^v, \alpha_n \rangle \xi_{n-1} \xi_n =
0.
\end{align*}
Also, the following relations hold \cite{8}:
\begin{align*}
&\psi_* (\xi_K) = 0\ \ \hbox{for length}\ \ (w_K) \neq |K|,\\[.2pc]
&\psi_* (\xi_K) = [X_{w_K w_0}]\ \ \hbox{if length}\ \ (w_K) =
|K|,
\end{align*}
where
\begin{equation*}
w_K = \prod\limits_{i\in K} s_{\alpha_i}, \quad K \subset [1,n].
\end{equation*}

In view of p.~58 of \cite{8}, the class of the tangent bundle of
$Z$ in $K(Z)$ can be explicitly described as
\begin{equation*}
[T_Z] = \sum_{r=1}^n \left[ \bigotimes\limits_{i=1}^n {\mathcal
O}\;(Z_i)^{\langle \alpha_i^v, \alpha_r \rangle} \right].
\end{equation*}
Accordingly, the Todd character of $T_Z$ is
\begin{equation}
td(Z) = \prod\limits_{r =1}^n \frac{\sum_{i=1}^r \langle
\alpha_i^v, \alpha_r \rangle \xi_i}{1 - \exp \left( - \sum_{i=1}^r
\langle \alpha_i^v, \alpha_r \rangle \xi_i \right)}.
\end{equation}

From now on $K = [1,\dots,k], k \leq n$ is fixed. From the exact
sequence,
\begin{equation*}
0 \rightarrow T_{Z_K}\rightarrow T_Z |Z_K \rightarrow N_{Z_K \mid
Z} \rightarrow 0,
\end{equation*}
the isomorphism $N_{Z_K\mid Z} \cong {\mathcal O}\;(Z_1)|Z_K
\oplus \cdots \oplus {\mathcal O}\;(Z_k)|Z_K$ and from the
Riemann--Roch--Grothendieck formula applied to the inclusion
$i_K\hbox{\rm :}\ Z_K \rightarrow Z$ one obtains
\begin{equation}
\ch({\mathcal O}_{Z_K}) = (-1)^k \prod\limits_{i=1}^k
\frac{\xi_i^2}{1-\e^{\xi_i}}.
\end{equation}

Notice that for every element $w\in W$ there exists a
decomposition of $w_0$ into simple reflections, such that
$w_{[1,\dots,k]} = ww_0$, where $k = \hbox{length}(w)$ (see
\cite{8}).

On the other hand, for $K = [1,\dots,k]$, one knows that \cite{8}
\begin{equation*}
\psi_*({{\mathcal O}}_{Z_K}) = {{\mathcal O}}_{X_w}\quad \hbox{\rm
and}\quad R^q \psi_* ({\mathcal O}_{Z_K}) = 0, \quad q \geq 1.
\end{equation*}

Finally, Riemann--Roch--Grothendieck formula yields
\begin{equation}
\ch({\mathcal O}_{X_w}) \e^{-\rho} = \psi_* (\ch({\mathcal
O}_{Z_K}) td(Z)).
\end{equation}

In conclusion, by substituting (A1) and (A2) into (A3) and by
knowing the action of $\psi_*$ on the cohomology, one gets an
algorithm of computing $\ch({\mathcal O}_{X_w})$ in terms of the
Cartan integers $\langle \alpha_i^v, \alpha_j \rangle$, a reduced
decomposition of $w$ and the generators $[X_{w'}], w'\in W$ of
$H^*(G/B, \Q)$. However, as mentioned before in \S5, some further
simplifications of this algorithm are worth having.

\begin{rema}{\rm It is likely that the recent paper by Kostant and
Kumar \cite{15} can lead to a better understanding of the
integrality conditions in \S5.}
\end{rema}

\section*{Acknowledgements}

The present article was completed in 1990, a few months before the
sudden and tragic death of the first author. The second author was
partially supported by a grant from the National Science
Foundation. The main results were announced in C~B\v{a}nic\v{a}
and M~Putinar, Fibr\'{e}s vectoriels complexes de rang stable sur
less vari\'{e}t\'{e}s complexes compactes, {\it C.~R.~Acad. Sci.
Paris.~t.~Serie~I} {\bf 314} (1992) 829--832. The second author
thanks his colleagues at the Indian Statistical Institute,
Bangalore and the Tata Institute, Mumbai for their interest in
ressurrecting this work.

\end{document}